\centerline{A FORMULA FOR THE k-th COVARIANT DERIVATIVE}
\medskip
\centerline{Kostadin Tren\v cevski}
\medskip
\centerline{Institute of Mathematics, Ss Cyril and Methodius University in Skopje,}
\centerline{Arhimedova 3, 1000 Skopje, Macedonia, e-mail: kostadin.trencevski@gmail.com}
\medskip
{\bf Abstract.} The aim of the present paper is to give a formula for  the  $k$-th  covariant 
derivative of tensor field along a given curve. In order to do that, first 
the symbols $P^{i<k>}_{j}$ and $Q^{i<k>}_{j}$ which depend on the 
Christoffel symbols are introduced. Some properties of them are also given. 
The main result is given by (3.1) and further it is generalized for $k\in R$. 
\medskip 
\par 
{\bf 1. Introduction of the symbols $P^{i<k>}_{j}$ and $Q^{i<k>}_{j}$.}
\par
Let x(s) be a given curve on a differentiable manifold $M_{n}$ which is 
endowed with affine connection $\Gamma $. The $k$-th covariant derivative
$\nabla ^{k}_{\dot x}$ along the curve $x(s)$ will be denoted by 
$\nabla ^{k}$. We shall deal further in a local coordinate system. 
\par 
The next three equalities 
$$P^{i<0>}_{j}=\delta ^{i}_{j}, \qquad P^{i<1>}_{j}=\Gamma ^{i}_{jl}
{dx^{l}\over ds},$$
$$P^{i<k+1>}_{j}  = {dP^{i<k>}_{j}\over ds} + P^{l<k>}_{j}P^{i<1>}_{l}, 
\qquad (k\in \{ 0,1,2,\cdots \}), \eqno{(1.1)}$$
determine all symbols $P^{i<k>}_{j}$ $(i,j\in \{1,\cdots ,n\})$, uniquely. 
Analogously we introduce the symbols $Q^{i<k>}_{j}$ by the following 
three equalities 
$$Q^{i<0>}_{j}=\delta ^{i}_{j}, \qquad Q^{i<1>}_{j}=-\Gamma ^{i}_{jl}
{dx^{l}\over ds},$$
$$Q^{i<k+1>}_{j} = {dQ^{i<k>}_{j}\over ds} + Q^{i<k>}_{l}Q^{l<1>}_{j}, 
\qquad (k\in \{ 0,1,2,\cdots \}), \eqno{(1.2)}$$
\medskip
\par 
{\bf 2. Some properties for the symbols $P^{i<k>}_{j}$ and $Q^{i<k>}_{j}$.}
\par 
(i) {\it The  transformation  law  for  the  components 
$P^{i<k>}_{j}$ from 
one coordinate system to another is given by the following formula
$$P^{j<l>}_{r} = \sum _{p=0}^{l}P'^{i<p>}_{q}
{\partial x^{j}\over \partial x'^{i}} 
\Bigl ({\partial x'^{q}\over \partial x^{r}} \Bigr )^{(l-p)} C^{p}_{l},
\eqno{(2.1)}$$
\noindent where $(l-p)$ denotes $(l-p)$-th derivative  by $s$.} 
\par 
{\it Proof.} We will prove (2.1) by induction of the index $l$. It is 
satisfied for $l=0$. Assuming that (2.1) holds for $l$, then by 
differentiating of (2.1) we obtain 
$$P^{j<l+1>}_{r} - P^{k<l>}_{r}P^{j<1>}_{k} = 
\sum _{p=0}^{l}
(P'^{i<p+1>}_{q}-P'^{a<p>}_{q}P'^{i<1>}_{a})
{\partial x^{j}\over \partial x'^{i}} 
\Bigl ({\partial x'^{q}\over \partial x^{r}} \Bigr )^{(l-p)} C^{p}_{l} $$
$$
+ \sum _{p=0}^{l}
P'^{i<p>}_{q}\Bigl [{d \over ds}{\partial x^{j}\over 
\partial x'^{i}}\Bigr ]
\Bigl ({\partial x'^{q}\over \partial x^{r}} \Bigr )^{(l-p)} C^{p}_{l} +
\sum _{p=0}^{l}
P'^{i<p>}_{q}{\partial x^{j}\over \partial x'^{i}}
\Bigl ({\partial x'^{q}\over \partial x^{r}} \Bigr )^{(l+1-p)} C^{p}_{l}, $$
$$
P^{j<l+1>}_{r} = 
\sum _{p=1}^{l+1}
P'^{i<p>}_{q}
{\partial x^{j}\over \partial x'^{i}} 
\Bigl ({\partial x'^{q}\over \partial x^{r}} \Bigr )^{(l+1-p)} C^{p-1}_{l}-
\sum _{p=0}^{l}
P'^{a<p>}_{q}P'^{i<1>}_{a}
{\partial x^{j}\over \partial x'^{i}}
\Bigl ({\partial x'^{q}\over \partial x^{r}} \Bigr )^{(l-p)} C^{p}_{l} $$
$$
+ \sum _{p=0}^{l}
P'^{i<p>}_{q}\Bigl [{d\over ds}{\partial x^{j}\over 
\partial x'^{i}}\Bigr ]
\Bigl ({\partial x'^{q}\over \partial x^{r}} \Bigr )^{(l-p)} C^{p}_{l} +
\sum _{p=0}^{l}
P'^{i<p>}_{q}{\partial x^{j}\over \partial x'^{i}}
\Bigl ({\partial x'^{q}\over \partial x^{r}} \Bigr )^{(l+1-p)} C^{p}_{l} $$
$$ +\Bigl [\sum_{p=0}^{l}P'^{c<p>}_{q}
{\partial x^{k}\over \partial x'^{c}}
\Bigl ({\partial x'^{q}\over \partial x^{r}} \Bigr )^{(l-p)} C^{p}_{l}\Bigr ]
\Bigl [ {\partial x^{j}\over \partial x'^{q}}\Bigl ({d\over ds}
{\partial x'^{q}\over \partial x'^{k}}\Bigr ) + P'^{i<1>}_{a}
{\partial x^{j}\over \partial x'^{i}}
{\partial x'^{a}\over \partial x^{k}}\Bigr ].$$
We notice that 
$$ P'^{i<p>}_{q}\Bigl ({d\over ds}
{\partial x^{j}\over \partial x'^{i}} \Bigr )
\Bigl ({\partial x'^{q}\over \partial x^{r}} \Bigr )^{(l-p)} C^{p}_{l} +
P'^{c<p>}_{q}{\partial x^{k}\over \partial x'^{c}}
\Bigl ({\partial x'^{q}\over \partial x^{r}} \Bigr )^{(l-p)} C^{p}_{l} 
{\partial x^{i}\over \partial x'^{q}}
\Bigl ({d\over ds}{\partial x'^{q}\over \partial x^{k}} \Bigr )  $$
$$ =
P'^{c<p>}_{q}
\Bigl ({\partial x'^{q}\over \partial x^{r}} \Bigr )^{(l-p)} C^{p}_{b} 
{\partial x^{k}\over \partial x'^{c}}
\Bigl [{\partial x'^{i}\over \partial x^{k}}
\Bigl ({d\over ds}{\partial x^{j}\over \partial x'^{i}} \Bigr ) +
{\partial x^{j}\over \partial x'^{i}}
\Bigl ({d\over ds}{\partial x'^{i}\over \partial x^{k}} \Bigr )\Bigr ] =0.$$
Further we obtain 
$$
P^{j<l+1>}_{r} = \sum _{p=1}^{l+1}P'^{i<p>}_{q}
{\partial x^{j}\over \partial x'^{i}}
\Bigl ({\partial x'^{q}\over \partial x^{r}}\Bigr )^{(l+1-p)}C^{p-1}_{l}+ $$
$$
+\sum _{p=0}^{l}P'^{i<p>}_{q}
{\partial x^{j}\over \partial x'^{i}}
\Bigl ({\partial x'^{q}\over \partial x^{r}}\Bigr )^{(l+1-p)}C^{p}_{l} =
\sum _{p=0}^{l+1}P'^{i<p>}_{q}
{\partial x^{j}\over \partial x'^{i}}
\Bigl ({\partial x'^{q}\over \partial x^{r}}\Bigr )^{(l+1-p)}C^{p}_{l+1}. $$
\par 
Analogously to the previous proof, one can prove the following 
property. 
\par 
(i) {\it The transformation law for the components 
$Q^{i<k>}_{j}$ from 
one coordinate system to another is given by the following formula}
$$Q^{j<l>}_{r} = \sum _{p=0}^{l}Q'^{a<p>}_{b}
{\partial x'^{b}\over \partial x^{r}} 
\Bigl ({\partial x^{j}\over \partial x'^{a}} \Bigr )^{(l-p)} C^{p}_{l}.
\eqno{(2.2)}$$
\par 
(iii) {\it The $k$-th derivative of $P^{i<l>}_{r}$ can be expressed in the 
following form }
$$(P^{j<l>}_{r})^{(k)} = \sum _{i=0}^{k}P^{a<l+k-i>}_{r}Q^{j<i>}_{a}C^{i}_{k}
. \eqno{(2.3)}$$
\par 
{\it Proof.} We will prove (2.3) by induction of $k$. It is satisfied for 
$k=0$. If it holds for the numbers $\{ 1,\cdots ,k\}$, then by 
differentiation of (2.3) we obtain 
$$(P^{j<l>}_{r})^{(k+1)} = \sum _{i=0}^{k} [P_{r}^{a<l+1+k-i>} - 
P_{r}^{\lambda <l+k-i>} P_{\lambda }^{a<1>}] Q^{j<i>}_{a} C^{i}_{k}$$
$$ + \sum _{i=0}^{k} P_{r}^{a<l+k-i>}[Q^{j<i+1>}_{a}- Q^{j<i>}_{\lambda }
Q^{\lambda <1>}_{a}] C^{i}_{k}. $$
Using that $P^{a<1>}_{\lambda }=-Q^{a<1>}_{\lambda }$, we obtain 
$$ (P_{r}^{j<l>})^{(k+1)} = \sum _{i=0}^{k}P_{r}^{a<l+k+1-i>}Q^{j<i>}_{a}
C^{i}_{k} + \sum _{i=1}^{k+1}P_{r}^{a<l+k+1-i>}Q^{j<i>}_{a}C^{i-1}_{k} = 
\sum _{i=0}^{k+1}P_{r}^{a<l+(k+1)-i>}Q^{j<i>}_{a}C^{i}_{k+1}. $$

Analogously to (iii), it holds 

(iv) 
$$(P^{j<l>}_{r})^{(k)} = \sum _{i=0}^{k}Q^{j<l+k-i>}_{a}P^{a<i>}_{r}C^{i}_{k}
. \eqno{(2.4)}$$

(v) {\it The following formula holds}
$$ \sum _{r=0}^{k}P_{u}^{j<k-r>} Q_{t}^{u<r>} C^{r}_{k} = 
\left \{ \matrix{ \delta ^{j}_{t}& if & k=0\cr
0& if & k\ge 1\cr }\right. .\eqno{(2.5)}$$

{\it Proof.} We will prove (2.5) by induction of the index $k$. It is 
trivially satisfied for $k=0$ and for $k=1$. Assume that (2.5) holds for 
the number $k\ge 1$. Then by differentiation of (2.5) we obtain 
$$ \sum _{r=0}^{k}(P_{u}^{j<k-r+1>} - P_{u}^{\lambda <k-r>}
P_{\lambda }^{j<1>})Q_{t}^{u<r>}C^{r}_{k} + 
\sum _{r=0}^{k}P_{u}^{j<k-r>} (Q_{t}^{u<r+1>}-Q_{\lambda }^{u<r>}
Q^{\lambda <1>}_{t})C^{r}_{k} = 0.$$
If we replace $r-1$ instead of $r$ in the sum 
$\sum\limits _{r=0}^{k}P^{j<k-r>}_{u}
Q^{u<r+1>}_{t} $ and use that $P^{i<1>}_{j}=-Q^{i<1>}_{j}$, we obtain that 
(2.5) holds for $k+1$ also. 
\medskip 
{\bf 3. Formula for the $k$-th covariant derivative.}

{\bf Theorem 1.} {\it The $k$-th covariant derivative of a tensor field 
$A\in C^{k}$ of type $(p,q)$ can be expressed in the following form}
$$\nabla ^{k} A^{j_{1}\cdots j_{p}}_{i_{1}\cdots i_{q}} = \sum _{a=0}^{k} 
\sum _{m_{1}+\cdots +m_{p}+l_{1}+\cdots +l_{q}=a} 
P_{u_{1}}^{j_{1}<m_{1}>} \cdots P_{u_{p}}^{j_{p}<m_{p}>}\times $$
$$ \times Q_{i_{1}}^{v_{1}<l_{1}>} \cdots Q_{i_{q}}^{v_{q}<l_{q}>} 
{d^{k-a}A^{u_{1}\cdots u_{p}}_{v_{1}\cdots v_{q}}\over ds^{k-a}}
{k!\over (k-a)!m_{1}!\cdots m_{p}!l_{1}! \cdots l_{q}!}. \eqno{(3.1)}$$
{\it where $m_{1},\cdots ,m_{p},l_{1},\cdots ,l_{q}$ take values in the set}
$\{0,1,2,\cdots \}$. 

{\it Remark.} We convenient to write $\sum m$ instead of 
$m_{1}+\cdots +m_{p}$ and $\sum l$ instead of $l_{1}+\cdots +l_{q}$. If we 
use that $(-l)!=\pm \infty $ for $l\in \{1,2,3,...\}$, then (3.1) can be 
written in the form 
$$\nabla ^{k} A^{j_{1}\cdots j_{p}}_{i_{1}\cdots i_{q}} = 
P_{u_{1}}^{j_{1}<m_{1}>} \cdots P_{u_{p}}^{j_{p}<m_{p}>}
Q_{i_{1}}^{v_{1}<l_{1}>} \cdots Q_{i_{q}}^{v_{q}<l_{q}>} \times $$
$$\times {d^{k-\Sigma m-\Sigma l}
A^{u_{1}\cdots u_{p}}_{v_{1}\cdots v_{q}}\over ds^{k-\Sigma m-\Sigma l}}
{k!\over (k-\Sigma m-\Sigma l)!
m_{1}!\cdots m_{p}!l_{1}! \cdots l_{q}!}, \eqno{(3.2)}$$
where it is understood that the symbols of summation over the indices 
$m_{1},\cdots ,m_{p},l_{1},\cdots ,l_{q}$ in the set $\{0,1,\cdots \}$ 
are omitted and the symbols of summation over the indices 
$u_{1},\cdots ,u_{p},v_{1},\cdots ,v_{q}$ in the set $\{1,\cdots ,n\}$ 
are also omitted. 

{\it Proof.} One can easily verify that $\nabla ^{0}$ is the identity 
operator and $\nabla ^{1}$ is the covariant differentiation along  the  curve 
$x(s)$. In order to prove the theorem it is sufficient to prove that
$$ \nabla (\nabla ^{k}A) = \nabla ^{k+1}A.\eqno{(3.3)}$$
Using (1.1), (1.2) and the definition of $\nabla $, we obtain 
$$\nabla (\nabla ^{k}A^{j_{1}\cdots j_{p}}_{i_{1}\cdots i_{q}}) = \Bigl \{
P^{j_{1}<m_{1}+1>}_{u_{1}}\cdots P^{j_{p}<m_{p}>}_{u_{p}}
Q^{v_{1}<l_{1}>}_{i_{1}}\cdots Q^{v_{q}<l_{q}>}_{i_{q}}+$$
$$+ \cdots +P^{j_{1}<m_{1}>}_{u_{1}}\cdots P^{j_{p}<m_{p}+1>}_{u_{p}}
Q^{v_{1}<l_{1}>}_{i_{1}}\cdots Q^{v_{q}<l_{q}>}_{i_{q}}+
P^{j_{1}<m_{1}>}_{u_{1}}\cdots P^{j_{p}<m_{p}>}_{u_{p}}
Q^{v_{1}<l_{1}+1>}_{i_{1}}\cdots Q^{v_{q}<l_{q}>}_{i_{q}}+$$
$$+ \cdots +P^{j_{1}<m_{1}>}_{u_{1}}\cdots P^{j_{p}<m_{p}>}_{u_{p}}
Q^{v_{1}<l_{1}>}_{i_{1}}\cdots Q^{v_{q}<l_{q}+1>}_{i_{q}}\Bigr \}\times $$
$$ \times {d^{k-\Sigma m-\Sigma l}
A^{u_{1}\cdots u_{p}}_{v_{1}\cdots v_{q}}\over ds^{k-\Sigma m-\Sigma l}}
{k!\over (k-\Sigma m-\Sigma l)!
m_{1}!\cdots m_{p}!l_{1}! \cdots l_{q}!} + $$
$$ + P^{j_{1}<m_{1}>}_{u_{1}}\cdots P^{j_{p}<m_{p}>}_{u_{p}}
Q^{v_{1}<l_{1}>}_{i_{1}}\cdots Q^{v_{q}<l_{q}>}_{i_{q}}\times $$
$$\times {d^{k+1-\Sigma m-\Sigma l}
A^{u_{1}\cdots u_{p}}_{v_{1}\cdots v_{q}}\over ds^{k+1-\Sigma m-\Sigma l}}
{k!\over (k-\Sigma m-\Sigma l)!
m_{1}!\cdots m_{p}!l_{1}! \cdots l_{q}!}. $$
If we replace $m_{1}$ instead of $m_{1}+1$ in the first summand, 
$m_{2}$ instead of $m_{2}+1$ in the second summand and so on, after summing
up we obtain (3.3). 
\medskip 
{\bf 4. Fractional covariant derivatives.}

Analogously to the formula (3.1) now we can introduce an operator 
$\nabla ^{\alpha }$ for $\alpha \in R$ by the following formula 
$$\nabla^{\alpha }A^{j_{1}\cdots j_{p}}_{i_{1}\cdots i_{q}} = 
\sum _{a=0}^{\infty } 
\sum _{m_{1}+\cdots +m_{p}+l_{1}+\cdots +l_{q}=a} 
P_{u_{1}}^{j_{1}<m_{1}>} \cdots P_{u_{p}}^{j_{p}<m_{p}>}\times $$
$$ \times Q_{i_{1}}^{v_{1}<l_{1}>} \cdots Q_{i_{q}}^{v_{q}<l_{q}>} 
{d^{\alpha -a}A^{u_{1}\cdots u_{p}}_{v_{1}\cdots v_{q}}\over 
ds^{\alpha -a}}
{\alpha (\alpha -1)\cdots (\alpha -a+1)
\over m_{1}!\cdots m_{p}!l_{1}! \cdots l_{q}!}. \eqno{(4.1)}$$
where $m_{1},\cdots ,m_{p},l_{1},\cdots ,l_{q}$ take values in the set
$\{0,1,2,\cdots \}$. If $a=0$ in (4.1), then we define 
$\alpha (\alpha -1)\cdots (\alpha -a+1)=1$. If the Christoffel symbols 
are identically equal to zero, then 
$$\nabla^{\alpha }A^{j_{1}\cdots j_{p}}_{i_{1}\cdots i_{q}} = 
d^{\alpha }A^{j_{1}\cdots j_{p}}_{i_{1}\cdots i_{q}}/ds^{\alpha }$$
and so (4.1) is a generalization of the ordinary fractional derivative [3].
In [1] it is given another formula for fractional covariant derivatives. 

We notice that the right side of (4.1) is a series. If that series diverges
then we should apply a method for summing up divergent series. 

It is obvious that $\nabla ^{\alpha }$ coincides with the $\alpha $-th 
covariant derivative if $\alpha \in \{0,1,2,\cdots \}$. In order to prove 
some other properties for the operator $\nabla ^{\alpha }$ we will use the 
following two properties for the ordinary fractional derivatives 
$${d^{\alpha }\over ds^{\alpha }}\circ {d^{\beta }\over ds^{\beta }} = 
{d^{\alpha +\beta }\over ds^{\alpha +\beta }},\eqno{(4.2)}$$
$${d^{\alpha }\over ds^{\alpha }}(f_{1}\cdots f_{k}) = 
\sum _{l_{1}=0}^{\infty }\sum _{l_{2}=0}^{\infty } \cdots 
\sum _{l_{k-1}=0}^{\infty }
{d^{l_{1}}f_{1}\over ds^{l_{1}}} {d^{l_{2}}f_{2}\over ds^{l_{2}}}\cdots 
{d^{l_{k-1}}f_{k-1}\over ds^{l_{k-1}}} \times $$
$$\times 
{d^{\alpha -l_{1}-\cdots -l_{k-1}}f_{k}\over 
ds^{\alpha -l_{1}-\cdots -l_{k-1}}}
{\alpha (\alpha -1)\cdots (\alpha -l_{1}-l_{2}-\cdots -l_{k-1}+1)\over 
l_{1}!l_{2}!\cdots l_{k-1}!}\eqno{(4.3)}$$
and we shall assume everywhere that the ordinary fractional derivatives 
are defined. We shall also assume in the next proofs that that the 
corresponding series are convergent and we convenient to omit the symbols 
of summation. 

(i)  $\nabla ^{\alpha }$ {\it maps the tensors of type (p,q) into 
tensors of the same type. }

{\it Proof.} We should prove the transformation formula 
$$ \nabla ^{\alpha }A^{j_{1}\cdots j_{p}}_{i_{1}\cdots i_{q}} = 
{\partial x^{j_{1}}\over \partial x'^{r_{1}}} \cdots 
{\partial x^{j_{p}}\over \partial x'^{r_{p}}}
{\partial x'^{a_{1}}\over \partial x^{i_{1}}} \cdots 
{\partial x'^{a_{q}}\over \partial x^{i_{q}}} 
\nabla ^{\alpha }A'^{r_{1}\cdots r_{p}}_{a_{1}\cdots a_{q}}. \eqno{(4.4)}$$
Using the formulas (2.1), (2.2) and (4.3), we obtain 
$$ \nabla ^{\alpha }A^{j_{1}\cdots j_{p}}_{i_{1}\cdots i_{q}} = 
P'^{a_{1}<c_{1}>}_{b_{1}} {\partial x^{j_{1}}\over \partial x'^{a_{1}}} 
\Bigl ({\partial x'^{b_{1}}\over \partial x^{u_{1}}}
\Bigr )^{(k_{1}-c_{1})} C^{c_{1}}_{k_{1}}\cdots 
P'^{a_{p}<c_{p}>}_{b_{p}} {\partial x^{j_{p}}\over \partial x'^{a_{p}}} 
\Bigl ({\partial x'^{b_{p}}\over \partial x^{u_{p}}}
\Bigr )^{(k_{p}-c_{p})} C^{c_{p}}_{k_{p}}\times $$
$$\times 
Q'^{m_{1}<n_{1}>}_{d_{1}} {\partial x'^{d_{1}}\over \partial x^{i_{1}}} 
\Bigl ({\partial x^{v_{1}}\over \partial x'^{m_{1}}}
\Bigr )^{(l_{1}-n_{1})} C^{n_{1}}_{l_{1}}\cdots 
Q'^{m_{q}<n_{q}>}_{d_{q}} {\partial x'^{d_{q}}\over \partial x^{i_{q}}} 
\Bigl ({\partial x^{v_{q}}\over \partial x'^{m_{q}}}
\Bigr )^{(l_{q}-n_{q})} C^{n_{q}}_{l_{q}} \cdots \times $$
$$\times 
\Bigl ({\partial x^{u_{1}}\over \partial x'^{r_{1}}}\Bigr )^{(t_{1})}\cdots 
\Bigl ({\partial x^{u_{p}}\over \partial x'^{r_{p}}}\Bigr )^{(t_{p})}
\Bigl ({\partial x'^{s_{1}}\over \partial x^{v_{1}}}\Bigr )^{(w_{1})}\cdots
\Bigl ({\partial x'^{s_{q}}\over \partial x^{v_{q}}}\Bigr )^{(w_{q})}\times
$$
$$\times {d^{\alpha -\Sigma k-\Sigma l-\Sigma t-\Sigma w}\over 
ds^{\alpha -\Sigma k-\Sigma l-\Sigma t-\Sigma w}}
A'^{r_{1}\cdots r_{p}}_{s_{1}\cdots s_{q}}
{(\alpha -\Sigma k-\Sigma l)\cdots (\alpha -\Sigma k-\Sigma l
-\Sigma t-\Sigma w +1)\over t_{1}!\cdots t_{p}!w_{1}! \cdots w_{q}!}\times $$
$$\times {\alpha (\alpha -1)\cdots (\alpha -\Sigma k -\Sigma l +1)\over 
k_{1}!\cdots k_{p}!l_{1}!\cdots l_{q}!}.\eqno{(4.5)}$$
We notice that for a fixed value of the index $c_{1}$, the last term can be 
written in the following form
$$\sum _{a=0}^{\infty } \sum _{k_{1}-c_{1}+t_{1}=a} 
\Bigl ({\partial x'^{b_{1}}\over \partial x^{u_{1}}}\Bigr )^{(k_{1}-c)}
{1\over (k_{1}-c_{1})!}
\Bigl ({\partial x^{u_{1}}\over \partial x'^{r_{1}}}\Bigr )^{(t_{1})} M(a)$$
such that for fixed value of $k_1 +t_1=a+c_1 $ the term $M(a)$ do not depend 
on $k_1$ and $t_1$. So that term is equal to 
$$\sum _{a=0}^{\infty } 
\Bigl ({\partial x'^{b_{1}}\over \partial x'^{r_{1}}}\Bigr )^{(a)}
{1\over a!}M(a) = \delta ^{b_{1}}_{r_{1}} M(0).$$
So in the right side of the equality (4.5) we should substitute 
$k_1 =c_1 $and $t_1  =0$. Analogously we should substitute there 
$k_2 =c_2 ,t_2 =0, \cdots , k_p = c_p , t_p =0, l_1 =n_1, w_1 =0,\cdots ,
l_q =n_q, w_q =0$. Then (4.4) can easily be obtained. 

(ii) $\nabla ^{\alpha }$ {\it commutes with the contraction.}

{\it Proof.} Let us assume that $j_p =i_q =t$ in (4.1). Using the property (2.3)
for $l=0$, we obtain 
$$\sum _{m_{p}=0}^{\infty }\sum _{l_{q}=0}^{\infty }
P^{t<m_{p}>}_{u_{p}}Q^{v_{q}<l_{p}>}_{t} {1\over m_{p}!}{1\over l_{q}!} = 
\sum _{a=0}^{\infty }\sum _{m_{p}+l_{q}=a}
P^{t<m_{p}>}_{u_{p}}Q^{v_{q}<l_{p}>}_{t} {1\over m_{p}!}{1\over l_{q}!} = 
\sum _{a=0}^{\infty }{1\over a!}(\delta ^{v_{q}}_{u_{p}})^{(a)} = 
\delta ^{v_{q}}_{u_{p}}.$$
Substituting  that in the right side of (4.1). we obtain the proof of this
property. 

(iii) 
$$\nabla ^{\beta }
\Bigl (\nabla ^{\alpha }A^{j_{1}\cdots j_{p}}_{i_{1}\cdots i_{q}}\Bigr )
=\nabla ^{\alpha +\beta }A^{j_{1}\cdots j_{p}}_{i_{1}\cdots i_{q}}.
\eqno{(4.6)}$$

{\it Proof.} Using the definition (4.1) twice and the formulas (4.2), (4.3),
(2.3) and (2.4), one obtains 
$$\nabla ^{\beta }
\Bigl (\nabla ^{\alpha }A^{j_{1}\cdots j_{p}}_{i_{1}\cdots i_{q}}\Bigr ) =
P^{j_{1}<k_{1}>}_{u_{1}}\cdots P^{j_{p}<k_{p}>}_{u_{p}}
Q^{v_{1}<l_{1}>}_{i_{1}}\cdots Q^{v_{q}<l_{q}>}_{i_{q}} \times $$
$$\times 
P^{d_{1}<a_{1}+c_{1}>}_{w_{1}}Q^{u_{1}<r_{1}-c_{1}>}_{d_{1}}C^{c_{1}}_{r_{1}}
\cdots 
P^{d_{p}<a_{p}+c_{p}>}_{w_{p}}Q^{u_{p}<r_{p}-c_{p}>}_{d_{p}}C^{c_{p}}_{r_{p}}
\times $$
$$\times 
Q^{t_{1}<b_{1}+n_{1}>}_{m_{1}}P^{m_{1}<s_{1}-n_{1}>}_{v_{1}}C^{n_{1}}_{s_{1}}
\cdots 
Q^{t_{q}<b_{q}+n_{q}>}_{m_{q}}P^{m_{q}<s_{q}-n_{q}>}_{v_{q}}C^{n_{q}}_{s_{q}}
\times $$
$$\times 
{d^{(\alpha -\Sigma a-\Sigma b)+(\beta -\Sigma k-\Sigma b)-\Sigma r-\Sigma s}
\over ds^{\alpha +\beta -\Sigma a-\Sigma b-\Sigma k-\Sigma l-\Sigma r-
\Sigma s}} A^{w_{1}\cdots w_{p}}_{t_{1}\cdots t_{q}} 
{\beta (\beta -1)\cdots (\beta -\Sigma k-\Sigma l-\Sigma r-\Sigma s +1)
\over k_1 !\cdots k_p !l_1 !\cdots l_q !}\times $$ 
$$\times {\alpha (\alpha -1)\cdots (\alpha -\Sigma a-\Sigma b +1)\over 
a_1 !\cdots a_p !b_1 !\cdots b_q !}
{1\over r_1 !\cdots r_p !s_1 !\cdots s_q !}. $$
Similarly as in the proof of (i) we assume that $c_1 $ is fixed. The last 
term can be separated into sums for which $k_1 +r_1 -c_1 =a(=const.)$ 
where $a$ changes from $0$ to infinity. From (2.5) we obtain 
$$  
P^{j_{1}<k_{1}>}_{a_{1}}Q^{u_{1}<r_{1}-c_{1}>}_{d_{1}}
{1\over k_{1}!(r_{1}-c_{1})!} = 
\left \{ \matrix{
\delta ^{j_{1}}_{d_{1}} & if\; r_{1}=c_{1}\; and \; k_{1}=0\cr
0& otherwise \cr }\right. .$$
Using that fact in the above equality, we should substitute there 
$r_1 =c_1 $ and $k_1 =0$. Analogously we should substitute there 
$r_2 =c_2 , k_2 =0, \cdots ,r_p =c_p , k_p =0, s_1 =n_1 ,
l_1 =0, \cdots ,s_q =n_q, l_q =0$. After these substitutions it is easy 
to see that if we prove the following formula 
$$ \sum _{a_1 +c_1 =k_1 } \cdots \sum _{a_p +c_p =k_p }
\sum _{b_1 +n_1 =l_1 }\cdots \sum _{b_q +n_q =l_q }
{\beta (\beta -1)\cdots (\beta -\Sigma c-\Sigma n +1)\over 
c_1 !\cdots c_p ! n_1 ! \cdots n_q !}\times $$
$$\times {\alpha (\alpha -1)\cdots (\alpha -\Sigma a-\Sigma b +1)\over 
a_1 !\cdots a_p ! b_1 ! \cdots b_q !} = 
{(\alpha +\beta )(\alpha +\beta -1)\cdots 
(\alpha +\beta -\Sigma k-\Sigma l +1)\over 
k_1 !\cdots k_p ! l_1 ! \cdots l_q !}, \eqno{(4.7)}$$
where $a_{1},\cdots ,a_{p}, c_{1},\cdots ,c_{p}, b_{1},\cdots ,b_{q},
n_{1},\cdots ,n_{q}\in \{0,1,2,\cdots \}$, then the proof of (4.6) will be 
finished. The formula (4.7) can be regarded as equality between two 
polynomials of $\alpha $ and $\beta $. So it is sufficient to prove (4.7) 
if $\alpha $ and $\beta $ are positive integers. Suppose that 
$f_{1},\cdots ,f_{p}, g_{1},\cdots ,g_{q}, h\in C^{\alpha +\beta }$ are 
arbitrary functions of $t$. The coefficients in front of 
$$f_{1}^{(k_{1})}\cdots f_{p}^{(k_{p})}g_{1}^{(l_{1})}\cdots g_{q}^{(l_{q})}
h^{(\alpha +\beta -k_{1}-\cdots -k_{p}-l_{1}-\cdots -l_{q})}$$
from the left and the right side of the equality 
$$[(f_1 \cdots f_p g_1 \cdots g_q h)^{(\alpha )}]^{(\beta )} = 
(f_1 \cdots f_p g_1 \cdots g_q h)^{(\alpha +\beta )}$$
should be equal, and as a consequence we obtain (4.7). 
\medskip 

{\bf 5. Applications.}

The formula (3.1) gives a simple method of calculation of 
the $k$-th covariant derivative along a given curve $x(s)$. Further, the 
formula (4.1) for $\alpha =-1$ gives a solution of covariant tensor integral. 
Our solution reduces to calculate all of the symbols $P^{i<k>}_{j}$ and 
$Q^{i<k>}_{j}$ $(i,j\in \{1,\cdots ,n\}, k\in \{0,1,\cdots \})$ and to 
calculate the successive integrals of the field $A$, which is easier than 
to solve a system of differential equations [1]. 

If we consider the following system of linear differential equations
$$dY^{i}/ds + \sum _{j=1}^{n}f^{i}_{j}Y^{j} = g^{i}, \quad i=1,\cdots ,n,
\eqno{(5.1)}$$
where $f^{i}_{j}$, $g^{i}$ are functions of $s$, then by putting 
$P^{i<1>}_{j}=f^{i}_{j}$, we obtain the following tensor equation
$\nabla ^{1}Y^{i}=g^{i}$ and its solution is $Y^{i}=\nabla ^{-1}g^{i}$. 

Similarly, if we consider the following system of linear differential 
equations 
$$d^{2}Y^{i}/ds^{2} + 2\sum _{j=1}^{n}f^{i}_{j}dY^{j}/ds + 
\sum _{r=1}^{n}[df^{i}_{r}/ds + \sum _{l=1}^{n}f^{i}_{l}f^{l}_{r}]Y^{r} = 
g^{i}, \quad i=1,\cdots ,n,\eqno{(5.2)}$$
where $f^{i}_{j}$, $g^{i}$ are functions of $s$, then by putting 
$P^{i<1>}_{j}=f^{i}_{j}$, we obtain the following tensorial equation
$\nabla ^{2}Y^{i}=g^{i}$ and its solution is $Y^{i}=\nabla ^{-2}g^{i}$. 
\medskip
\centerline{REFERENCES}
\item {[1]} W.Fabian, Tensor Integrals, Proc. of the Edinburgh Math. Soc., 
(2) 10, part IV, 145-151. 
\item {[2]} S.Kobayashi, K.Nomizu, Foundations of Differential Geometry, 
Vol. 1, Interscience Publishers, New York 1963. 
\item {[3]} S.G.Samko, A.A.Kilbas, O.I.Marichev, Integrals and fractional 
derivatives and some of their applications, Minsk 1987, (in Russian). 
\end